\input amssym.def
\input amssym.tex

\hfuzz=15pt

\def\ii{{\frak I}}
\def\ia{{\frak a}}
\def\ib{{\frak b}}
\def\im{{\frak m}}
\def\reg{\hbox{\rm reg}}

\def\indeg{\hbox{\rm indeg}}
\def\proj{\hbox{\rm Proj}}
\def\Z{{\cal Z}}
\def\C{{\cal C}}

\def\lra{\longrightarrow}

\def\max{\hbox{\rm max}}

\def\fini{{$\quad\quad\square$}}

\def\reg{\hbox{\rm reg}}

\def\indeg{\hbox{\rm indeg}}

\def\proj{\hbox{\rm Proj}}

\def\tor{\hbox{\rm Tor}}

\font\nrm=cmcsc10 at10pt

\font\eightrm=cmr10 scaled 900
\font\eightit=cmti10 scaled 900
\font\bigrm=cmb10 scaled 1600
\font\mbigrm=cmb10 scaled 1200

\bigskip\bigskip\bigskip\bigskip

{\bigrm Sur la r\'egularit\'e de Castelnuovo-Mumford des id\'eaux,}

{\bigrm  en dimension 2}

\bigskip
\medskip

{\bf Marc CHARDIN et Amadou Lamine FALL}
\bigskip
\bigskip\bigskip
\bigskip

{\bf R\'esum\'e.} {\eightrm }  Nous montrons une borne pour la r\'egularit\'e de Castelnuovo-Mumford d'un id\'eal homog\`ene $I$ d'un anneau de polyn\^omes $A$ en termes du nombre de
variables et des degr\'es des g\'en\'erateurs dans le cas o\`u la dimension de $A/I$ est au plus deux.  Cette borne am\'eliore celle obtenue par Caviglia et Sbarra dans [CS]. Puis, en s'inspirant de l'article [CD], nous construisons \`a partir de  familles  de courbes monomiales
des id\'eaux homog\`enes  ayant une r\'egularit\'e proche des bornes fournies pr\'ec\'edemment. 
\bigskip
\centerline{\bf On the Castelnuovo-Mumford regularity of ideals, in dimension 2}
\medskip

{\bf Abstract.} {\eightit }We give a  bound on the
Castelnuovo-Mumford regularity of a homogeneous ideals $I$, in a
polynomial ring $A$, in terms of number of variables and the degrees of
generators, when the dimension of $A/I$ is at most two.  This bound improves the one obtained by
Caviglia and Sbarra in [CS]. In the continuation of the examples constructed in [CD], we use families of monomial curves to construct homogeneous ideals showing that these bounds are quite sharp.

\bigskip

\bigskip\bigskip
{\mbigrm Abridged English Version}\medskip

In this Note we give a refinement of a bound due to Caviglia and Sbarra [CS] on the Castelnuovo-Mumford regularity of homogeneous ideals $I$, in a polynomial ring $A$ over a field, when $\dim (A/I)=2$. If $A$ is of dimension $m+2$, we prove :\medskip

{\nrm Theorem 1.1.} {\it
Let $I$ be a homogeneous $A$-ideal  of codimension $m$ generated in degrees $d_{1}\geq d_{2}\ldots\geq d_{s}$ with $s>m$. Let 
$\Z$ be the zero-dimensionnal scheme defined by a general linear section of  $X:=\proj (A/I)$, and $i_{\Z}=\indeg
(I_{\Z})$ be the smallest degree of a hypersurface containing $\Z$. Then,
$$
\reg(I)\leq (d_{1}d_{2}\cdots d_{m}-\deg (\Z )+1
)(d_{1}+d_{2}+\ldots+d_{m+1}-m-i_{\Z})+i_{\Z}.
$$
}
\medskip
This and previous results in the cases  of dimensions zero and  one gives
\medskip
{\nrm Corollary 1.3.}   {\it If $I\subset A$ is a $A$-ideal generated in degrees at most $d$, then
\smallskip
{\rm (i)}  $\reg(A/I)\leq (m+2)(d-1)$ if $\dim (A/I)\leq 1$,
\smallskip
{\rm (ii)} $\reg(A/I)\leq (m+1)d^{m}(d-1)$ if $\dim (A/I)=2$.}
\medskip

The bound in (i) is known to be essentially optimal. We provide examples showing that this is also the case for the bound in (ii), at least regarding the power of $d$ involved.

These examples are a continuation of the ones in [CD, 2.3] and [CD, 2.5]. They use the monomial curve 
$\C_{m,n}\subset {\bf  P}^{m+2}_{k}$ parametrized on the affine chart $X_{0}=1$ by $(1:t:t^{n^{m}}:t^{n^{m-1}(n+1)}:\ldots:t^{(n+1)^{m}})$ and  its projection $\C'_{m,n}$ on $\{ X_{m+2}=0\} \simeq {\bf  P}^{m+1}_{k}$. 

We provide complete intersection ideals $I_{m,n}$ and $I'_{m,n}$ defining respectively the union of $\C_{m,n}$ and two fat lines and the union of $\C'_{m,n}$ and two fat lines. We then use an estimate on the regularity of the union of the two fat lines, to produce almost complete intersection ideals $\ii_{m,n}$ and $\ii'_{m,n}$ satisfying the following propositions :

{\nrm Proposition 2.2.} {\it The ideal $\ii'_{m,n}\subset k[X_{0},\ldots ,X_{m+2}]$ is generated by $m$ forms of degree $n+1$, one form of degree $2^{m-1}+1$ and one form of degree $mn+2^{m-1}$. One has $\dim(A/\ii'_{m,n})=2$ and
$$
\reg(\ii'_{m,n})\geq n^{m}+mn+2^{m-1}-1.
$$
}

{\nrm Proposition 3.2.} {\it 
The ideal $\ii_{m,n}\subset k[X_{0},\ldots ,X_{m+1}]$ is generated by $m-1$ forms of degree $n+1$, one form of degree $2^{m-2}+n$ and one form of degree $mn+2^{m-2}-1$. One has $\dim(A/\ii_{m,n})=2$ and
$$
\reg(\ii_{m,n})\geq n^{m}+mn+2^{m-2}-2.
$$
}

In the setting of Proposition 3.2 the bound given by Theorem 1.1 shows that
$$
\reg(\ii_{m,n})\leq 2m^{2}n(n+1)^{m-2}(n+2^{m-2})^{2}.
$$

The exponent of $n$ in this expression is $m+1$ which is one more than its actual value. \medskip

{\bf 1. Un raffinemment de la borne de Caviglia et Sbarra, en dimension 2}\medskip

Dans toute cette section, $A=k[X_{0},\ldots ,X_{m+1}]$ d\'esigne un anneau de polyn\^omes en $m+2$ variables sur un corps $k$ et $\im :=(X_{0},\ldots ,X_{m+1})$. Si $M$ est un $A$-module gradu\'e de type fini, $\indeg (M)$ d\'esigne le minimum des degr\'es des  \'el\'ement de $M$ si $M\not= 0$. La r\'egularit\'e de Castelnuovo-Mumford de $M$ est d\'efinie par
$$
\reg (M)=\max_{i}\{ a_{i}(M)+i\} =\max_{i}\{ b_{i}(M)-i\} .
$$
avec les notations $a_{i}(M):=\max\{ \mu\ |\ H^{i}_{\im}(M)_{\mu}\not= 0\}$,
si $H^{i}_{\im}(M)\not= 0$ et $a_{i}(M):=-\infty$ sinon et $b_{i}(M):=\max \{ \mu \ |\ \tor_{i}^{R}(M,k)_{\mu}\not= 0\} $
si $\tor_{i}^{R}(M,k)\not= 0$ et $b_{i}(M):=-\infty$ sinon. Si $I$ est un id\'eal homog\`ene de $A$, on pose $I^{sat}:=\cup_{j>0}(I:\im^{j})$ de telle sorte que $I^{sat}/I=H^{0}_{\im}(A/I)$. Si $\Z =\proj (A/I)$ est un sous-sch\'ema de ${\bf P}^{m+1}_{k}$, on pose $I_{\Z}:=I^{sat}$. 
\medskip

En utilisant la m\'ethode de Caviglia et Sbarrra, et un r\'esultat de [Ch] en dimension 1, nous obtenons le r\'esultat plus pr\'ecis suivant :\medskip

{\nrm Th\'eor\`eme 1.1.} {\it
Soit $I$ un ideal de $A$ de codimension $m$, engendr{\'e}
en degr{\'e}s $d_{1}\geq d_{2}\ldots\geq d_{s}$ avec $s>m$. Soit
$\Z$ le groupe de points d{\'e}fini par une section
g{\'e}n{\'e}rale de $X:=\proj (A/I)$, et $i_{\Z}:=\indeg
(I_{\Z})$. Alors,
$$
\reg(I)\leq (d_{1}d_{2}\cdots d_{m}-\deg (\Z )+1
)(d_{1}+d_{2}+\ldots+d_{m+1}-m-i_{\Z})+i_{\Z}.
$$
}

Notons que $\deg (\Z )=\deg (X)$ et que si $\C$ est la composante de
dimension 1 de $X$ et $I_{\C}$ son ideal de d\'efinition, $d_{s}\geq \indeg (I_{\C})\geq i_{\Z}\geq 1$.\medskip

Pour la preuve, nous aurons besoin du lemme suivant qui am\'eliore la borne de [CS, 2.3] \medskip

{\nrm Lemme 1.2.} {\it Soit $I$ un id{\'e}al
homog{\`e}ne de $A$, tel que $I\not= I:\im$, $l$ une forme lin{\'e}aire telle que $I:(l)/I$ soit de longueur finie et $q$ le plus petit entier tel que $I:l^{\infty}=I:l^{q}$.
 Alors on a $q\leq a_{0}(A/I)-\indeg(I^{sat})+1\leq \reg (I)-\indeg(I^{sat})$.
 }

{\it Preuve du Lemme 1.2.} La seconde in{\'e}galit{\'e} est claire car $a_{0}(A/I)+1\leq \reg (A/I)+1=\reg (I)$. Soit $z\in
I:l^{\infty}=I^{sat}$, et posons $K:=a_{0}(A/I)+1-\indeg (I^{sat})$.
Il s'agit de montrer que $zl^{K}\in I$. On peut
supposer $z$ homog{\`e}ne et non nul, ce qui entra{\^\i}ne que
$\deg z\geq \indeg(I^{sat})$. On a $zl^{K}\in (I^{sat})_{K+\deg
(z)}$ et $K+\deg (z)\geq a_{0}(A/I)+1$. Or $I_{\mu}=(I^{sat})_{\mu}$
pour $\mu \geq a_{0}(A/I)+1$. Il en d{\'e}coule que
$(I^{sat})_{K+\deg (z)}=I_{K+\deg (z)}$, d'o{\`u}  $zl^{K}\in I$.
\fini\medskip

{\it Preuve du Th\'eor\`eme 1.1.} Soit $l_{m+1},l_{m+2}$  des
formes lin{\'e}aires telles que $(I:(l_{m+2}))/I$ et ${A}/(I+(l_{m+1},l_{m+2}))$ soient de longueurs finies. On suppose de plus que $l_{m+2}$ est assez
g{\'e}n{\'e}rale, de telle sorte que posant $\Z :=\proj
(A/I+(l_{m+2}))$ on soit dans les conditions du th{\'e}or{\`e}me.
(Notons que l'on a $I_{\Z}=(I+(l_{m+2}))^{sat}$, par
d{\'e}finition).

Appliquons alors les r{\'e}sultats de la preuve de ([CS, 2.4]) avec
 $c=m$. Comme $I$ est engendre en degr{\'e}s $\leq
d_{1}$, en posant  $R:=\reg(I+(l_{m+2}))$  on a, par [CS, (2.1)],
$$
\reg(I)\leq
\max\{d_{1},R\}+\lambda\left( {I:l_{m+2}}\over{I}
\right)
$$
o{\`u} $\lambda(M)$ d{\'e}signe la longueur d'un $A$-module $M$.
De plus, d'apr{\`e}s [CS, (2.2)] et le Lemme 1.2, on a
$$
\lambda\left( {{I:l_{m+2}}\over{I}} \right) \leq
\lambda\left( {{I_{\Z}+(l_{m+1})}\over{I+(l_{m+1},l_{m+2})}} \right)
(R-i_{\Z}).
$$

D'autre part,
$$
\lambda\left(
{I_{\Z}+(l_{m+1})}\over{I+(l_{m+1},l_{m+2})}\right)
=\lambda\left( {A}\over{I+(l_{m+1},l_{m+2})}\right)
-\lambda\left( {A}\over{I_{\Z}+(l_{m+1})}\right)
\leq d_{1}\cdots d_{m}-\deg (\Z ).
$$
car $I+(l_{m+1},l_{m+2})$ contient une intersection compl{\`e}te de degr{\'e}s
$d_{1},\ldots ,d_{m},1,1$.

On a donc au total :
$$
\reg(I)\leq \max\{d_{1},R\}+(d_{1}\ldots
d_{m}-\deg (\Z ))(R-i_{\Z}).
$$

D'apr{\`e}s [Ch, 3.3], $R\leq d_{1}+\cdots +d_{m+1}-m$, d'o{\`u}
la borne annonc{\'e}e. \fini\medskip

{\nrm Corollaire 1.3.}   {\it Soit $I\subset A=k[X_{0},\ldots,X_{m+1}]$ un id{\'e}al engendr{\'e} en
degr{\'e}s au plus $d$. Alors,\smallskip

{\rm (i)}  $\reg(A/I)\leq (m+2)(d-1)$
si $\dim (A/I)\leq 1$, \smallskip

{\rm (ii)} $\reg(A/I)\leq
(m+2)d^{m}(d-1)$ si $\dim (A/I)=2$
}
\medskip
{\it Preuve.} Le cas (i) d\'ecoule de [Ch, 3.3], et le (ii) du Th\'eor\`eme 1.1 une fois not\'e que $\deg (\Z )\geq 1$ et $i_{\Z}\geq 1$.\fini

\bigskip

{\bf 2. Premier exemple d'id\'eaux}\medskip

Soit $m\geq 1$, $A'=k[X_{0},\ldots,X_{m+2}]$, $\ib'_{m,n}$ l'id{\'e}al de la courbe
monomiale
$\C'_{m,n}\subset {\bf P}^{m+2}$ param\'etr\'ee sur la carte affine $X_{0}=1$ par $(1:t:t^{n^{m}}:t^{n^{m-1}(n+1)}:\ldots:t^{(n+1)^{m}})$.

Les polyn\^omes $F_{i}=X_{i}^{n+1}-X_{0}X_{i+1}^{n}$ pour $2\leq i\leq m+1$ sont dans $\ib'_{m,n}$.

D\'efinissons par r\'ecurrence sur $m$, $f_{1}(X_{2},X_{3}):={{X_{3}}\over{X_{2}}}$ et $f_{m}(X_{2},\ldots ,X_{m+2}):=f_{m-1}({{X_{3}}\over{X_{2}}},\ldots, {{X_{m+2}}\over{X_{m+1}}})$. On note que si $m=1$, et $x=(x_{0}:\cdots :x_{3})\in \C'_{1,n}$ avec $x_{0}x_{1}\not= 0$, alors $f_{1}(x_{2},x_{3})={{x_{1}}\over{x_{0}}}$. D'autre part, si $m>1$ et $x=(1:x_{1}:\cdots :x_{m+2})\in \C'_{m,n}$ avec $x_{1}\not= 0$, alors $(1:x_{1}:{{x_{3}}\over{x_{2}}}:\cdots :{{x_{m+2}}\over{x_{m+1}}})\in \C'_{m-1,n}$; on en d\'eduit par r\'ecurrence sur $m$ que $f_{m}(x_{2},\ldots ,x_{m+2})=x_{1}$. De plus, on v\'erifie, par r\'ecurrence sur $m$, que
$$
f_{m}(X_{2},\ldots ,X_{m+2})^{(-1)^{m}}=\prod_{j=0}^{m}X_{j+2}^{(-1)^{j}{{m}\choose{j}}}={{P_{m}(X_{2},X_{4},\ldots )}\over {Q_{m}(X_{3},X_{5},\ldots )}}
$$
o\`u $P_{m}$ et $Q_{m}$ sont des polyn\^omes homog\`enes de degr\'e $2^{m-1}$. Il en d\'ecoule que le polyn\^ome d\'efini par $F'_{1}:=X_{0}P_{m}-X_{1}Q_{m}$ si $m$ est pair et $F'_{1}:=X_{1}P_{m}-X_{0}Q_{m}$ si $m$ est impair appartient \`a $\ib'_{m,n}$. Soit $I'_{m,n}:=(F'_{1},F_{2},\ldots ,F_{m+1})$. \medskip

{\nrm Lemme 2.1.} {\it On a $I'_{m,n}=\ib'_{m,n}\cap J'_{m,n}\cap K'_{m,n}$, o\`u $J'_{m,n}$ est un id\'eal $(X_{2},\ldots ,X_{m+2})$-primaire et
$K'_{m,n}$ est un id\'eal  $(X_{0},X_{2},\ldots ,X_{m+1})$-primaire. En particulier $I'_{m,n}$ est une intersection compl\`ete de codimension $m+1$.
}\medskip

{\it Preuve.} Soit $J':=(X_{2},\ldots ,X_{m+2})$ et $K':=(X_{0},X_{2},\ldots ,X_{m+1})$. Le sch\'ema  $\Z'_{m,n}:=\proj (A'/I'_{m,n})$ contient les droites
$D:=\proj (A'/J')$ et
$D':=\proj (A'/K')$. De plus
$\Z'_{m,n}\cap \{X_{0}=0\}$ est support{\'e} sur $D'$.

Montrons  que sur la carte affine $X_{0}=1$, $\Z'_{m,n}$ est l'union d'un sch\'ema support\'e sur $D$ et de
$\C'_{m,n}$.  On note que si $X_{2}\not= 0$, alors $X_{i}\not= 0$ pour
tout $i$.  Lorsque $m=1$, on a $X_{2}^{n+1}=X_{3}^{n}=\left( X_{1}X_{2}\right) ^{n}$ d'o\`u $X_{2}=X_{1}^{n}$ puis $X_{3}=X_{1}X_{2}=X_{1}^{n+1}$.
Si $m>1$, on pose $Y_{1}:=X_{1}$ et $Y_{i}:={{X_{i+1}}\over{X_{i}}}$ pour $2\leq i\leq m+1$ et on remarque que $X_{i}=Y_{i}^{n}$ pour $2\leq i\leq m+1$, d'o\`u $Y_{i}^{n+1}=Y_{i+1}^{n}$ pour $2\leq i\leq m$. Par d\'efinition $f_{m-1}(Y_{2},\ldots ,Y_{m+1})=f_{m}(X_{2},\ldots ,X_{m+2})$. Par r\'ecurrence sur $m$ on en d\'eduit que $Y_{i}=Y_{1}^{n^{m-i+1}(n+1)^{i-2}}$, d'o\`u $X_{i}=Y_{i}^{n}=Y_{1}^{n^{m-i+2}(n+1)^{i-2}}=X_{1}^{n^{m-i+2}(n+1)^{i-2}}$ pour $2\leq i\leq m+1$.
De plus $X_{m+2}=X_{m+1}Y_{m+1}=X_{1}^{n(n+1)^{m-1}}Y_{1}^{(n+1)^{m-1}}=X_{1}^{(n+1)^{m}}$. Ceci montre, comme annonc\'e, que $\Z'_{m,n}$ co\"\i ncide avec $\C'_{m,n}$ hors de $D\cup D'$.

$\Z'_{m,n}$ est donc de codimension $m+1$. Ainsi
$(F'_{1},F_{2},\ldots F_{m+1})$ est n\'ecessairement une suite r\'eguli\`ere dans $A'$ et $I'_{m,n}$ est purement de codimension $m+1$. Le lemme en d\'ecoule.\fini\medskip

L'id\'eal $\ib'_{m,n}$ \'etant radical, on a $\reg (J'_{m,n}\cap K'_{m,n})<\reg (I'_{m,n})=mn+2^{m-1}+1$ par [CU, 4.2]. En particulier, $\ia'_{m,n}:=J'_{m,n}\cap K'_{m,n}$ est engendr\'e en degr\'es au plus $mn+2^{m-1}$. Il existe donc $F'\in \ia'_{m,n}-\ib'_{m,n}$ homog\`ene de degr\'e $d':=mn+2^{m-1}$. On pose $\ii'_{m,n}:=I'_{m,n}+(F')$.\medskip

{\nrm Proposition 2.2.} {\it L'id{\'e}al $\ii'_{m,n}$ est engendr{\'e} par $m$ polyn{\^o}mes
de degr{\'e} $n+1$, un polyn{\^o}me de degr{\'e} $2^{m-1}+1$ et un
polyn{\^o}me de degr{\'e} $mn+2^{m-1}$. On a $\dim
(A'/\ii'_{m,n})=2$ et
$$
\reg(\ii'_{m,n})\geq n^{m}+mn+2^{m-1}-1.
$$
}

{\it Preuve.} On a une suite
exacte $0\lra A'/\ib'_{m,n}[-d'] \buildrel{\times F'}\over{\lra} A'/I'_{m,n}\lra A'/\ii'_{m,n} \lra 0$ d'o{\`u}
on d{\'e}duit que $H_{\im}^{0}(A'/\ii'_{m,n})\simeq
H_{\im}^{1}(A'/\ib'_{m,n})[-d']$ car $A'/I'_{m,n}$ est Cohen-Macaulay de dimension 2. En particulier
$$
a_{0}(A'/\ii'_{m,n})=a_{1}(A'/\ib'_{m,n})+d'.
$$
On note que
$X_{1}^{n^{m}}-X_{0}^{n^{m}-1}X_{2}$ est un g{\'e}n{\'e}rateur
minimal de $\ib'_{m,n}$. Donc $\reg(\ib'_{m,n})\geq n^{m}$. Comme
d'autre part pour $(m,n)$ distinct de $(2,2)$, $(3,2)$ et $(4,2)$,
$$
a_{2}(A'/\ib'_{m,n})+2<a_{2}(A'/I'_{m,n})+2=\reg (A'/I'_{m,n})=mn+2^{m-1}\leq n^{m}-1
$$
on a $a_{1}(A'/\ib'_{m,n})+1=\reg (A'/\ib'_{m,n})\geq n^{m}-1$ dans ces cas. Lorsque $m=2,3,4$ et $n=2$ on v\'erifie directement en utilisant le logiciel Macaulay 2 [GS] que $a_{1}(A'/\ib'_{m,n})+1=\reg (A'/\ib'_{m,n})= n^{m}-1$. La proposition en d\'ecoule.\fini
\medskip

\medskip

{\bf 3. Second exemple d'id\'eaux}\medskip

Soit $m\geq 2$, $A=k[X_{0},\ldots,X_{m+1}]$ et $\ib_{m,n}=\ib'_{m,n}\cap A$ l'id{\'e}al de la courbe
monomiale
$\C_{m,n}\subset {\bf P}^{m+1}$ parametr\'ee sur la carte affine $X_{0}=1$ par $(1:t:t^{n^{m}}:t^{n^{m-1}(n+1)}:\ldots:t^{n(n+1)^{m-1}})$, projection de la courbe $\C'_{m,n}$ sur l'hyperplan $X_{m+2}=0$.

Les polyn\^omes $F_{i}=X_{i}^{n+1}-X_{0}X_{i+1}^{n}$ pour $2\leq
i\leq m$ sont dans $\ib_{m,n}$. Soit
$F_{1}:=X_{1}^{n}P_{m-1}-X_{0}^{n}Q_{m-1}$ si $m$ est pair et $F_{1}:=X_{0}^{n}P_{m-1}-X_{1}^{n}Q_{m-1}$ si $m$ est impair. On remarque que
$F_{1}\in \ib_{m,n}$ et on pose $I_{m,n}:=(F_{1},\ldots ,F_{m})$.
\medskip

{\nrm Lemme 3.1.} {\it On a $I_{m,n}=\ib_{m,n}\cap J_{m,n}\cap K_{m,n}$, o\`u $J_{m,n}$ est un id\'eal $(X_{2},\ldots ,X_{m+1})$-primaire et
$K_{m,n}$ est un id\'eal  $(X_{0},X_{2},\ldots ,X_{m})$-primaire. En particulier $I_{m,n}$ est une intersection compl\`ete de codimension $m$.
}\medskip

{\it Preuve.} Soit $J:=(X_{2},\ldots ,X_{m+1})$ et
$K:=(X_{0},X_{2},\ldots ,X_{m})$. Le sch\'ema  $\Z_{m,n}:=\proj
(A/I_{m,n})$ contient les droites $\Delta :=\proj (A/J)$ et
$\Delta':=\proj (A/K)$, et  $\Z_{m,n}\cap \{X_{0}=0\}$ est
support{\'e} sur $\Delta'$.

Montrons  que sur la carte affine $X_{0}=1$, $\Z_{m,n}$ est
l'union d'un sch\'ema support\'e sur $\Delta$ et de $\C_{m,n}$.  Si
$X_{2}\not= 0$, alors $X_{i}\not= 0$ pour tout $i$. Lorsque $m=2$,
on a $X_{2}^{n+1}=X_{3}^{n}=\left( X_{1}^{n}X_{2}\right) ^{n}$
d'o\`u $X_{2}=X_{1}^{n^{2}}$ puis
$X_{3}=X_{1}^{n}X_{2}=X_{1}^{n(n+1)}$. Si $m>2$, on pose
$Y_{1}:=X_{1}$ et $Y_{i}:={{X_{i+1}}\over{X_{i}}}$ pour $2\leq
i\leq m$ et on remarque que $X_{i}=Y_{i}^{n}$ pour $2\leq i\leq m$
d'o\`u $Y_{i}^{n+1}=Y_{i+1}^{n}$ pour $2\leq i\leq m-1$. Par
d\'efinition $f_{m-2}(Y_{2},\ldots ,Y_{m})=f_{m-1}(X_{2},\ldots
,X_{m+1})$. Par r\'ecurrence sur $m$ on en d\'eduit que
$Y_{i}=Y_{1}^{n^{m-i+1}(n+1)^{i-2}}$, d'o\`u
$X_{i}=Y_{i}^{n}=Y_{1}^{n^{m-i+2}(n+1)^{i-2}}=X_{1}^{n^{m-i+2}(n+1)^{i-2}}$
pour $2\leq i\leq m$. De plus
$X_{m+1}=X_{m}Y_{m}=X_{1}^{n^{2}(n+1)^{m-2}}Y_{1}^{n(n+1)^{m-2}}=X_{1}^{n(n+1)^{m-1}}$.

$\Z_{m,n}$ est donc  support\'e sur la r\'eunion de $\C_{m,n}$,
$\Delta$ et $\Delta'$, et ainsi $(F_{1},\ldots ,F_{m})$ forme
n\'ecessairement une suite r\'eguli\`ere dans $A$. On conclue
comme pour le lemme 2.1 .\fini\medskip

 On a $\reg (J_{m,n}\cap
K_{m,n})<\reg (I_{m,n})=mn+2^{m-2}$, en particulier
$\ia_{m,n}:=J_{m,n}\cap K_{m,n}$ est engendr\'e en degr\'es au
plus $mn+2^{m-2}-1$. Il existe donc $F\in \ia_{m,n}-\ib_{m,n}$
homog\`ene de degr\'e $d:=mn+2^{m-2}-1$.\medskip
 En posant $\ii_{m,n}:=I_{m,n}+(F)$, on a \medskip

{\nrm Proposition 3.2.} {\it L'id{\'e}al $\ii_{m,n}$ est
engendr{\'e} par $m-1$ polyn{\^o}mes de degr{\'e} $n+1$, un
polyn{\^o}me de degr{\'e} $2^{m-2}+n$ et un polyn{\^o}me de
degr{\'e} $mn+2^{m-2}-1$. On a $\dim (A/\ii_{m,n})=2$ et
$$
\reg(\ii_{m,n})\geq n^{m}+mn+2^{m-2}-2.
$$
}
{\it Preuve.} On a une suite
exacte $0\lra A/\ib_{m,n}[-d] \buildrel{\times F}\over{\lra} A/I_{m,n}\lra A/\ii_{m,n} \lra 0$ d'o{\`u}
on d{\'e}duit que $a_{0}(A/\ii_{m,n})=a_{1}(A/\ib_{m,n})+d$ car $A/I_{m,n}$ est Cohen-Macaulay de dimension 2. 
On note que
$X_{1}^{n^{m}}-X_{0}^{n^{m}-1}X_{2}$ est un g{\'e}n{\'e}rateur
minimal de $\ib_{m,n}$. Donc $\reg(\ib_{m,n})\geq n^{m}$. Comme
d'autre part si $(m,n)$ est distinct de $(2,2)$ et $(3,2)$,
$$
a_{2}(A/\ib_{m,n})+2<a_{2}(A/I_{m,n})+2=\reg (A/I_{m,n})=mn+2^{m-2}\leq n^{m}-1
$$
on a $a_{1}(A/\ib_{m,n})+1=\reg (A/\ib_{m,n})\geq n^{m}-1$ dans ces cas. Lorsque $m=2$ ou $m=3$,  on a $a_{1}(A/\ib_{m,n})+1=\reg (A/\ib_{m,n})=n^{m}-1$ d'apr\`es [BCFH] si $m=2$ et  [CD, 2.6 (3)] si $m=3$. La proposition en d\'ecoule.\fini\medskip

{\nrm Remarque 3.3.} Lorsque $m=2$ ou $m=3$ on a plus pr\'ecis\'ement
$$
\reg(\ii_{m,n})=n^{m}+mn+2^{m-2}-2,
$$
car $\reg (\ib_{m,n})=a_{1}(A/\ib_{m,n})+2=n^{m}$ d'apr\`es [BCFH] si $m=2$ et  [CD, 2.6 (3)] si $m=3$. 
D'autre part, d'apr\`es [L, 5.5], $\reg (\ib_{m,n})\leq n^{m}+n(n+1)^{m-2}-1$ pour tout $m$, et ainsi
$$
\reg(\ii_{m,n})\leq n^{m}+n(n+1)^{m-2}+mn+2^{m-2}-3.
$$

Ceci doit \^etre compar\'e \`a la borne founie par le Th\'eor\`eme 1.1, qui donne, apr\`es quelques simplifications,
$$
\reg(\ii_{m,n})\leq 2m^{2}n(n+1)^{m-2}(n+2^{m-2})^{2}.
$$
L'exposant de $n$ dans cette expression diff\`ere donc de un par rapport \`a sa vraie valeur.

\bigskip\bigskip
\centerline{\bf R\'{e}f\'{e}rences bibliographiques}
\bigskip

\noindent [BS] D. Bayer, M. Stillman. A criterion for detecting
$m$-regularity.  {\it Invent. Math.}  {\bf 87} (1987), 1-11.\medskip

\noindent [BCFH] H. Bresinsky, F. Curtis, M.  Fiorentini, L. T. Hoa. On the structure of local cohomology modules for monomial curves in $P^{3}_{K}$. {\it Nagoya Math. J.}  {\bf 136}  (1994), 81--114.\medskip

\noindent [CS] G. Caviglia, E. Sbarra. Characteristic-free
bounds for the Castelnuovo-Mumford regularity, {\it pr\'epublication} math.AC/0310122.\medskip

\noindent [CD]  M. Chardin, C. D'Cruz. Castelnuovo-Mumford regularity: examples of curves and
surface. {\it J. Algebra}  {\bf 270} (2003), 347-360.\medskip

\noindent [CU]  M. Chardin,B. Ulrich. Liaison and the Castelnuovo-Mumford regularity.
 {\it Amer. J. Math.}  {\bf 124} (2002), 1103-1124.\medskip

\noindent [Ch] M. Chardin.  Regularity of ideals and their powers. {\it Pr\'epublication} 364. Institut de math\'e\-matiques de Jussieu, Mars 2004.
 \medskip

\noindent [GS] D. Grayson, M. Stillman. {\it Macaulay 2 software}, http://www.math.uiuc.edu/Macaulay2/. 
 \medskip
 
\noindent  [L] S. L'vovsky. On inflection points, monomial curves, and hypersurfaces containing projective curves.  {\it Math. Ann.}  {\bf 306}  (1996),  no. 4, 719--735.

\bigskip\bigskip\bigskip

\noindent M. C. : Institut de math\'ematiques de Jussieu,

\noindent CNRS et Universit\'e Paris VI,

\noindent 4, place Jussieu F-75005 Paris\smallskip

\noindent  A. L. F. : D\'epartement de math\'ematiques,

\noindent Facult\'e des Sciences, Universit\'e Cheikh Anta Diop,

\noindent Dakar, S\'en\'egal

\end